\documentclass[a4paper,11pt]{article}

\usepackage{graphicx}%
\usepackage{multirow}%
\usepackage{amsmath,amssymb,amsfonts}%
\usepackage{amsthm}%
\usepackage{mathrsfs}%
\usepackage[title]{appendix}%
\usepackage{xcolor}%
\usepackage{natbib}%
\usepackage{mathabx}
\usepackage[colorlinks=true, breaklinks, pdfstartview=FitV, linkcolor=green, citecolor=blue, urlcolor=gray,breaklinks]{hyperref}

\usepackage{times}

\usepackage{soul} 
\setstcolor{red}
\usepackage{cancel} 

\allowdisplaybreaks
\begin{document}

\title{Naturalization of Vinti's intermediary orbit---or how to get rid of small divisors in Brouwer's satellite solution}

\author{Martin Lara \\
Scientific Computation Research Institute (SCRIUR) University of La Rioja \\ \tt{mlara0@gmail.com}
 }


\date{\color{gray}Submitted  to Celest Mech Dyn Astron; August 5, 2026}

\maketitle

\begin{abstract}
A naturalization, in the sense of Deprit, of Vinti's dynamical system makes the construction of a first order analytical satellite perturbation solution without critical divisors possible. The process gives hints on how to recast Brouwer's classical solution in a form that is free from small divisors too. The functionality of classical perturbation solutions is in this way extended to the close vicinity of the critical inclination. The new analytical solution bears terms that may grow unbounded, but their negative effects do not disclose in time intervals in which this kind of solution may be useful in practice.
\end{abstract}

\section{Introduction}

Simple analytical solutions to the artificial satellite problem dating back to the beginning of the space era \citep{Sterne1958,Garfinkel1958,Vinti1959,Aksnes1965,CidLahulla1969}, still find application these days in spacecraft trajectory design or onboard orbit propagation \citep{GurfilLara2014CeMDA,Duffy2016,Wright2016,Senecal2023}. Beyond the practicable Keplerian orbits they take Earth's oblateness effects into account, and are gathered under the general denomination of \emph{main problem intermediaries}. Research on intermediary orbits continues today, either in the search for more significant models \citep{Deprit1981,HautesserresLaraCeMDA2017,BiriaRussell2018,Lara2019CMDA,BiriaRussell2020} or with the aim of refining computational algorithms for the efficient evaluation of existing intermediaries \citep{Vinti1998,Wiesel2015,Biria2023,Shepperd2024}. 
\par

Effective intermediary orbits replicate the general features of the mean, perturbed dynamics \citep{Aksnes1965}. In particular, they succeed in detecting two different regions of orbits with direct and retrograde motion of the perigee, respectively. The limit case in which the rotation of the perigee vanishes establishes a boundary of ``frozen orbits'' separating these two regions.
However, it is well known that minor effects of the main problem dynamics, or comparable small perturbations produced by the higher order harmonics of the gravitational potential, split the solution space into different regions of orbits with rotating and oscillating periapses \citep{PettyBreakwell1960,Struble1961,Kozai1961,Hagihara1961,Izsak1962}. Since these little effects are ignored by main problem intermediaries, the approximations provided by intermediary orbits do not capture the behavior of orbits with oscillating periapses---a flaw that is commonly attributed to their integrable character \citep{DepritMiller1988,CoffeyDepritDeprit1994,Lara2026Ke}.
\par

Simple analytical solutions exist too that model the Earth's shape more accurately than main problem intermediaries can do, and hence better approach the actual dynamics \citep{Garfinkel1959,Kozai1959,Brouwer1959}. But because these perturbation solutions are based on the complete reduction of the perturbation Hamiltonian to its secular terms, they neither cover the full range of solutions. Therefore, their application is constrained to the case of orbits with rotating line of apsides, which is the only instance in which the use of averaging operations is justified for the argument of the periapsis \citep{Tisserand1868}. Since both kinds of orbits---with oscillating an rotating perigees---exist about the Earth, separate analytical solutions are needed to tackle each distinct orbital regime \citep{Hori1960,Izsak1962,Garfinkel1966,Jupp1988}. 
\par

On the other hand, the small value of the Earth's oblateness makes that the stripes separating the regions with direct and retrograde rotation of the perigee, where the orbits with oscillating perigees may exist, are certainly narrow. As a consequence, the main efforts in the development of operational tools have been devoted to the analytical modeling of orbits with rotating perigees \citep{HootsRoehrich1980,CoffeyDeprit1982,CoffeyAlfriend1984,Solomon1991,ValladoCrawfordHujsakKelso2006,Fraysseetal2012,LaraSanJuanHautesserres2018}. Remark that this fact in no way denies the interest of satellite missions based on orbits with oscillating perigees \citep{Shapiro1996,Ulybyshev2009}.
\par

An additional trouble regarding analytical perturbation solutions obtained by total Hamiltonian reduction, which is not shared by intermediary orbits, is the appearance in the formal solution of divisors that vanish at the so-called \emph{critical inclinations} $i_\mathrm{c}=\frac{\pi}{2}\pm\arccos\sqrt{4/5}$. They ensue from the elimination of the long-period terms related to the perigee dynamics. Still, the singularity should not be much of worry because the region about the critical inclination is mostly filled with orbits with oscillating perigees, a case that, as it was already asserted, is not covered by these kinds of analytical solutions. But the appearance of critical divisors has undesired effects on orbits with rotating perigees too. To wit, the series representing the perturbation solution loose accuracy when approaching the critical inclination, even to the extent of turning their predictions into illusory \citep{Brouwer1959}. Protective limits for avoiding this sort of degradation in the convergence of the series are then stablished in practice. They roughly encompass inclinations in the interval ranging from 61 to 66 degree, or a corresponding interval of retrograde orbits \citep{CoffeyNealSegermanTravisano1995}. 
\par

The existence of elliptic orbits with oscillating perigees about the critical inclination is determined by the physical nature of the Earth's potential. Accordingly, it has been claimed that the problematic divisors ``cannot be removed by any method'' \citep{Hori1960}. The general acceptance of this premise turned the concerns on the evaluation of analytical perturbation solutions just into an issue of computational nature \citep{Solomon1991,HootsSchumacherGlover2004}. However, to our knowledge, Hori's (\citeyear{Hori1960}) categorical statement lacks supporting mathematical evidence. Thus, the question of whether the inability to remove the critical divisors from rotating-perigee perturbation solutions is intrinsic to the problem, as claimed, or, on the contrary, is instead related to the widely accepted philosophy of the perturbation methods used at the time remains open. Developing additional efforts to further elucidate this point seems, therefore, sufficiently justified.
\par

In this regard, it is worth noting that intermediary solutions are free from offending divisors. Then, finding a transformation that, up to some truncation order, turns the original Hamiltonian into the Hamiltonian of a main problem intermediary would prevent the trouble of the small divisors. More generally, without constraining to the case of intermediary orbits, the transformation would operate a preparatory \emph{simplification} on the original Hamiltonian \citep{Deprit1981,DepritFerrer1989,DepritMiller1989}. In this way the simplified Hamiltonian would act as a surrogate for the artificial satellite problem that could be approached by perturbations without the trouble of the critical divisors. Finding this kind of proxy might notably improve the accuracy of classical perturbation models in predicting the behavior of rotating-perigee orbits in the vicinity of the critical inclination.
\par

The desired transformation certainly exists. The inspiration for finding it lays on a recent demonstration regarding Vinti's dynamical system \citep{Vinti1959,Izsak1960}. Namely, that the models obtained from increasing truncations of Vinti's expanded Hamiltonian in spherical harmonics give rise to corresponding, successive, radical reductions of the size of the region of elliptic orbits with oscillating perigees \citep{Lara2026Ke}. Beyond usual applications of Vinti's analytical solution as a \emph{common} intermediary of the artificial satellite problem, such transformation would assign to Vinti's dynamical system the role of a \emph{natural} intermediary in the sense of \citet{Deprit1981}.
\par

In this naturalization process, an infinitesimal contact transformation is computed that turns a simplified gravitational model comprising only the perturbations due to the Earth's second and fourth zonal harmonics into the first significative truncation of Vinti's dynamical system. This kind of simplification does not affect formally the model, but simply removes one of the parameters by slightly tuning the value of the zonal harmonic coefficient of degree four. Namely, the transformation replaces $J_4$ by $-J_2^2$, where $J_n$ denotes a harmonic coefficient of degree $n$. This is precisely the condition identified by \citet{Brouwer1959} to cancel the terms with critical divisors out of his perturbation solution. Therefore, the new transformation is next applied to Brouwer's gravitational model. More precisely, as customary in practice, the contribution of the fifth zonal harmonic is ignored due to the smallness of the Earth's $J_5$ coefficient \citep{LaneCranford1969,HootsRoehrich1980,Hoots1982,ValladoCrawfordHujsakKelso2006}.
\par

It happens that, up to $\mathcal{O}(J_2^2)$, the removal of long-period terms related to the third zonal harmonic does not give rise to critical divisors. Therefore, the complete normalization of the simplified Hamiltonian in Brouwer's style leads to an analytical solution that is free from critical divisors in its periodic components. The secular terms of the new solution differ from those of Brouwer, but this is only through minor modifications of the coefficients defining the inclination polynomials.
\par

Regrettably, yet not surprisingly, the differences in the secular terms of both perturbation solutions have unavoidable consequences, and the proposed \emph{naturalization} of Vinti's intermediary involves secular terms, a fact that is certainly a weak point of the approach. However, at difference from the naturalization of usual common intermediaries of the main problem, for which the infinitesimal contact transformation is $\mathcal{O}(J_2)$, the fine tuning of $J_4$ is based on an $\mathcal{O}(J_2^2)$ transformation. As a consequence, the amplitudes of the unbounded terms needed to recover the osculating orbit grow at a moderate rate. Moreover, the transformation is further improved by ignoring the secular terms of the generating function yielding the naturalization of Vinti's dynamical model. In this way the terms that cause the small-divisor trouble are likewise removed, but still at the expense of dealing with unbounded terms. Luckily, they are now of the mixed secular-periodic type \citep{Poincare1892vII}, which, while being a concern in the development of classical perturbation theories \citep{Garfinkel1959,GarfinkelHoriAksnes1970}, enjoy a more benign character. A following perturbation approach  results in a new analytical solution that undergoes exactly the same secular variations as those of Brouwer's original solution.
\par

The performance of the new theory in the propagation of orbits close to the critical inclination is illustrated by comparison with corresponding reference orbits computed numerically for the same perturbation model. As expected, it shows radical improvements with respect to the classical approach, providing the awaited accuracy from this kind of perturbation solution. Far away enough from the critical inclination, the new theory provides similar accuracy to the one obtained with Brouwer's original approach in time intervals useful for common satellite applications, in which the mixed terms do not degrade appreciably the performance of the solution. Undeniably, the harmful effects of the unbounded terms will disclose eventually. However, this fact should not be a major worry in practical applications. In the end, the scope in which the accuracy provided by Brouwer's perturbation solution may be acceptable in practice, as it happens also to the one of Vinti's exact solution, is limited in time due to the insufficient modeling of the actual dynamics.
\par

\section{Vinti's dynamical system as a natural intermediary} \label{}

Vinti's (\citeyear{Vinti1959}) intermediary orbit is mainly celebrated for taking account of the whole contribution of the zonal harmonic of the second degree of the gravitational potential---the so-called ``main problem'' of artificial satellite theory \citep{Brouwer1959}.
\par

When Vinti's closed form Hamiltonian in oblate spheroidal coordinates is replaced by the usual expansion in spherical harmonics, one obtains what \cite{Izsak1960} called Vinti's dynamical system. That is,
\begin{equation} \label{Hoscu} 
\mathcal{H}=\sum_{m\ge0}\frac{J_2^m}{m!}\mathcal{H}_{m,0},
\end{equation}
where the oblateness coefficient $J_2$ plays the role of the small parameter. The integrable, zeroth-order term
\begin{equation} \label{Keppolar}
\mathcal{H}_{0,0}=\frac{1}{2}\left(R^2+\frac{\Theta^2}{r^2}\right)-\frac{\mu}{r}
\end{equation}
is the Keplerian in polar variables $(r,\theta,R,\Theta)$, standing for radius, argument of the latitude, radial velocity, and total angular momentum per unit of mass, respectively, whereas $\mu$ denotes the gravitational parameter. The perturbation terms are arranged in the form
\begin{equation} \label{Hm0}
\mathcal{H}_{m,0}=-\frac{\mu}{r}\frac{R_\Earth^{2m}}{r^{2m}}\sum_{j=0}^mb_{m,j}\cos2j\theta,
\end{equation}
where $R_\Earth$ denotes the Earth's mean equatorial radius, and $b_{m,j}\equiv{b}_{m,j}(i)$ are inclination polynomials. Referring to the canonical coordinates $i=\arccos(N/\Theta)$, where $N$ is the conjugate momenta to the longitude of the ascending node $\nu$.
\par

We only consider the first significant truncation of Eq.~(\ref{Hoscu}). Namely,
\begin{equation} \label{Vinti2}
\mathcal{H}=\mathcal{H}_{0,0}+J_2\mathcal{H}_{1,0}+\frac{1}{2}J_2^2\mathcal{H}_{2,0},
\end{equation}
in which
\begin{align} \label{H10}
\mathcal{H}_{1,0}=&-\frac{\mu}{r}\frac{R_\Earth^2}{r^2}\frac{1}{4}(2-3s^2+3s^2\cos2\theta), \\ \label{H20}
\mathcal{H}_{2,0}=&-\frac{\mu}{r}\frac{R_\Earth^4}{r^4}\frac{1}{32}\left[3(35s^4-40s^2+8)-20s^2(7s^2-6)\cos2\theta+35s^4\cos4\theta\right], 
\end{align}
where the abbreviation $s\equiv\sin{i}$ is used to shorten printed expressions.
\par

On the other hand, when the actual geopotential expansion is limited to contribution of the second and fourth zonal harmonics of the gravity potential, we have the perturbation Hamiltonian
\begin{equation} \label{geo2even}
\mathcal{K}=\mathcal{K}_{0,0}+J_2\mathcal{K}_{1,0}+\frac{1}{2}J_2^2\mathcal{K}_{2,0},
\end{equation}
in which $\mathcal{K}_{0,0}=\mathcal{H}_{0,0}$, $\mathcal{K}_{1,0}=\mathcal{H}_{1,0}$, and $\mathcal{K}_{2,0}=-\tilde{J}_4\mathcal{H}_{2,0}$, where we denote
\begin{equation} 
\tilde{J}_4\equiv\frac{J_4}{J_2^2}<0.
\end{equation}
\par

The slight difference between the perturbation Hamiltonians $\mathcal{K}$ and $\mathcal{H}$ has, nevertheless, radical consequences in the normalization of each perturbation Hamiltonian. Namely, the removal of long-period terms introduces critical divisors in the solution of $\mathcal{K}$, whereas it does not in the solution of $\mathcal{H}$; cf.~\citep{Brouwer1959}.

\subsection{Naturalization}

An infinitesimal contact transformation is now computed that, in the new variables, turns $\mathcal{K}$ into $\mathcal{H}$. Analogously to the naturalization of main-problem intermediaries advocated by \citet{Deprit1981}, it is derived from a generating function $\mathcal{W}$. The latter is the solution of the partial differential equation
\begin{equation} \label{natuhomo}
\{\mathcal{W},\mathcal{K}_{0,0}\}=\mathcal{K}-\mathcal{H},
\end{equation}
where $\mathcal{K}_{0,0}$ is the Keplerian term, and the curly brackets denote the Poisson brackets operator. When evaluated in polar coordinates,  Eq.~(\ref{natuhomo}) is readily solved by indefinite integration, to obtain
\begin{equation} \label{homonatu}
\mathcal{W}=\int\left(\mathcal{K}-\mathcal{H}\right)\frac{r^2}{\Theta}\mathrm{d}\theta.
\end{equation}
\par

Straightforward computations yield
\begin{equation} \label{Wnatu}
\mathcal{W}=-\Theta{J}_2^2\delta^2\frac{R_\Earth^4}{p^4}\left[(\mathcal{S}+\mathcal{M})\theta+\mathcal{P}\right],
\end{equation}
where $p=\Theta^2/\mu$,
\begin{equation}
\delta^2=-(1+\tilde{J}_4)\approx\frac{3}{8}
\end{equation}
for the Earth, and
\begin{align} \label{secular}
\mathcal{S}=&\; \frac{3}{128}(35 s^4-40 s^2+8)(3 \kappa ^2+3 \sigma ^2+2) \\ \label{natumixed}
\mathcal{M}=&-\frac{15}{64}s^2(7s^2-6) \left[\left(\kappa^2-\sigma ^2\right)\cos2\theta+2\kappa\sigma\sin2\theta\right], \\ \nonumber
\mathcal{P}=&\;\frac{1}{128}(35s^4-40s^2+8)(4\sigma^2+6\kappa^2+9\kappa+18)\sigma+\frac{1}{256}s^2(7s^2-6)\big\{2(16\sigma^2-8\kappa^2 \\ \nonumber
            & -15\kappa+40)\sigma\cos2\theta -[32\kappa(\kappa^2+2\sigma^2+5)+5(15\kappa^2+9\sigma^2+8)]\sin2\theta\big\} \\ \nonumber
            & +\frac{1}{512}[32\kappa(3\kappa^2-2\sigma^2+7)+35(7\kappa^2-\sigma^2+2)]s^4\sin4\theta \\ \label{natuperi}
            & +\frac{1}{128}(4\sigma^2-26\kappa^2-35\kappa-14)\sigma{s}^4\cos4\theta,
\end{align}
in which
\begin{equation} \label{skpolar}
\kappa=-1+\frac{r}{p}, \qquad
\sigma=\frac{Rp}{\Theta},
\end{equation}
denote the projections of the eccentricity vector in the orbital frame.
\par

Alternatively, the generating function of the naturalization can be written in terms of the canonical Delaunay variables. Thus,
\begin{align} \label{secular2}
\mathcal{S}=&\; \frac{3}{128}(35 s^4-40 s^2+8)(2+3e^2), \\ \label{natumixed2}
\mathcal{M}=&-\frac{15}{64}s^2(7s^2-6)e^2\cos(4f+2\omega), \\ \nonumber
\mathcal{P}=&\;\frac{1}{256}(35s^4-40s^2+8)[e^3\sin3f+9e^2\sin2f+9(e^2+4)e\sin{f}]
+\frac{1}{256}(7s^2-6) \\ \nonumber
&\times s^2\big\{2e^3\big[5\sin(-f+2\omega)-\sin(5f+2\omega)\big]-15e^2\sin(4f+2\omega)
-10(e^2+4)e \\ \nonumber
&\times\big[3\sin(f+2\omega)+\sin(3f+2\omega)\big]-20(3e^2+2)\sin(2f+2\omega)
\big\}
+\frac{1}{512}s^4\big\{ 5e^3  \\ \nonumber
& \times\big[7\sin(f+4\omega)+\sin(7f+4\omega)\big]+105e^2\sin(2f+4\omega)+35e^2\sin(6f+4\omega) \\\label{natuperi2}
            &+7(e^2+4)e\big[5\sin(3f+4\omega)+3\sin(5f+4\omega)\big]+35(3e^2+2)\sin(4f+4\omega) \big\},
\end{align}
where the eccentricity $e$, the argument of the periapsis $\omega$, and the true anomaly $f$ are functions of the Delaunay variables.
However, it is known that shorter transformation equations are commonly obtained in artificial satellite theories when formulated in polar variables \citep{Kozai1962,Izsak1963}. It happened to be the case in the current approach too, in spite of the compact expressions for $\mathcal{S}$ and $\mathcal{M}$ in Eqs.~(\ref{secular2}) and (\ref{natumixed2}), respectively.
\par

The generating function $\mathcal{W}$ in Eq.~(\ref{Wnatu}) is intentionally split into three parts to highlight the effects of the explicit appearance of $\theta$, which harm the ensuing infinitesimal contact transformation with secular terms as well as mixed secular-periodic terms of the form $\theta\sin2k\theta$ and $\theta\cos2k\theta$, $k=1,2$. While this flaw does not invalidate the transformation as such, it is definitely a weak point that constrains the validity of the natural intermediary to the short time intervals in which $\theta(\mathcal{S}+\mathcal{M})$ may remain below acceptable limits. 
\par

On the other hand, the fact that $\mathcal{W}$ is just $\mathcal{O}(J_2^2)$ certainly helps in delaying the degradation of the solution due to the effect of the unbounded terms. Additionally, it is worth noting that the correction $\Delta{a}=\{a,\mathcal{W}\}$ to the semimajor axis, namely
\begin{equation}
\{a,\mathcal{W}\}=-aJ_2^2\delta^2\frac{R_\Earth^4}{r^4}\frac{a}{r}\frac{1}{32}\left[3(35s^4-40s^2+8)-20(7s^2-6)s^2\cos2\theta+35s^4\cos4\theta\right],
\end{equation}
is not affected of secular or mixed terms, which is certainly very good news because $a$ is the more sensitive element regarding initialization of analytical or semi-analytical solutions to perturbed Kepler motions \citep{BreakwellVagners1970,Lara2020arxiv}. These two characteristics make that this kind of transformation can be useful in practice.
\par

Incidentally, secular terms appear too in the naturalization of Garfinkel's (\citeyear{Garfinkel1958}) main problem intermediary---refer to the term $W^\dag_g$ on p.~143 of \citep{Deprit1981}---yet constrained to the transformation of the argument of the latitude. Remarkably, they neither corrupt the infinitesimal contact transformation of the semimajor axis in that case, in spite of the transformation that naturalizes Garfinkel's intermediary is only $\mathcal{O}(J_2)$.
\par

It must be mentioned that the word ``naturalization'' is loosely used here, insofar as the truncation of Vinti's dynamical system in Eq.~(\ref{Vinti2}), not being integrable, is not an intermediary strictly speaking. However, it would be if the tail that has been cut to Vinti's dynamical system were added to both Eqs.~(\ref{Vinti2}) and (\ref{geo2even}). In spite of that, the word is used to highlight the parallelism of Deprit's procedure with the current approach. A generalization of the latter by means of Lie transforms, without constraint to a simple infinitesimal contact transformation, would more properly fit to the algebraic domain of canonical simplifications.
\par

\subsection{Mean dynamics}

The naturalization process turned the original,  $J_2$ and $J_4$ Hamiltonian in Eq.~(\ref{geo2even}) into Vinti's truncated Hamiltonian
\begin{equation} \label{Vinti21}
\mathcal{H}^\dag=\mathcal{H}_{0,0}+J_2\mathcal{H}_{1,0}+\frac{1}{2}J_2^2\mathcal{H}_{2,0},
\end{equation}
where $\mathcal{H}_{0,0}$ is the same given in Eq.~(\ref{Keppolar}), and $\mathcal{H}_{1,0}$ and $\mathcal{H}_{2,0}$ remain the same as those given in Eqs.~(\ref{H10}) and (\ref{H20}), respectively, yet all of them are, of course, now written in the new variables. In spite of the shortcomings associated to this surrogate due to the secular terms that corrupt the preliminary simplification, the exercise of computing a perturbation solution to Eq.~(\ref{Vinti21}) will be rewarding as discussed later at the beginning of section \ref{s:morsig}.
\par

A first-order analytical solution of Eq.~(\ref{Vinti21}) is then computed by Hamiltonian normalization in Brouwer's style. That is, the solution is made of two different sets of short- and long-period corrections, which are computed to the first order of $J_2$, and a set of secular terms, which are solved up to $\mathcal{O}(J_2^2)$. It is approached with the method of Lie transforms \citep{Deprit1969} taking $J_2$ as the small parameter. Rather than relying on the polar-nodal set of canonical variables, as it was previously done with the naturalization, the normalization is approached in Delaunay canonical variables, as customary. Namely, the mean anomaly $\ell=M$, the argument of the perigee $g=\omega$, and the longitude of the ascending node $h=\nu=\Omega$, for the coordinates, and $L=\sqrt{\mu{a}}$, $G=\Theta$, and $H=N$, for their respective conjugate momenta.

First, the short-period terms are removed from the truncation of Vinti's dynamical system by a Lie transformation
$T_1:(\ell,g,h,L,G,H;\epsilon)\rightarrow(\ell',g',h',L',G',H')$
that yields a Hamiltonian in mean elements. Next, a second Lie transformation
$
T_2:(\ell',g',h',L',G',H';\epsilon)\rightarrow(\ell'',g'',h'',L'',G'',H'')$,
removes the long-period terms from the mean dynamics, keeping only the secular terms. After truncation to $\mathcal{O}(J_2^2)$, the Hamiltonian is completely reduced to a function of the momenta only, and, therefore, it is trivially integrable.
\par

Up to $\mathcal{O}(J_2^2)$, the first step yields the new Hamiltonian
\begin{equation} \label{Vinti22}
\mathcal{H}^\ddag=\mathcal{H}_{0,0}+J_2\mathcal{H}_{0,1}+\frac{1}{2}J_2^2\mathcal{H}_{0,2},
\end{equation}
where the Keplerian remains unchanged, and
\begin{align} \label{H01}
\mathcal{H}_{0,1}=&\; \frac{\mu}{2a}\frac{R_\Earth^2}{p^2}\eta\frac{1}{2}\left(3s^2-2\right), 
\\ \nonumber
\mathcal{H}_{0,2}= &-\frac{\mu}{2a}\frac{R_\Earth^4}{p^4}\eta\frac{3}{16}\left[
5(21s^4-28s^2+8)+2(3s^2-2)^2\eta-2(25s^4-32s^2+8)\eta^2 \right. \\ \label{H02}
& \left. -4(5s^2-4)s^2e^2\cos2\omega \right],
\end{align}
where $\eta=G/L$, and $e=\sqrt{1-\eta^2}$; cf.~\citep{Lara2026Ke}.
\par

It must be kept in mind that, after Eq.~(\ref{Vinti22}) has been computed, the Hamiltonian flow corresponds to the mean dynamics, and hence all the symbols in Eqs.~(\ref{H01}) and (\ref{H02}) represent functions of the mean (prime) canonical variables. Additionally, the removal of short-period terms has turned the mean semimajor axis $a=L'^2/\mu$ into an integral of the mean dynamics. Together with the original integral $N=H$, they reduce Eq.~(\ref{Vinti22}) to a Hamiltonian of just one degree of freedom in $(g',G')$. Regarding the mean dynamics, $H$ it is customarily replaced by the constant inclination of the circular orbits $i_\mathrm{circular}=\arccos(H/\sqrt{\mu{a}})$.
\par

For comparison purposes, we also carry out an analogous removal of the short period terms from the $J_2$ and $J_4$ problem in Eq.~(\ref{geo2even}). It yields
\begin{equation} \label{Kprime}
\mathcal{K}'=\mathcal{K}_{0,0}+J_2\mathcal{K}_{0,1}+\frac{1}{2}J_2^2\mathcal{K}_{0,2},
\end{equation}
whose terms differ from those of Eq.~(\ref{Vinti22}) only in the term
\begin{equation} \label{K02}
\mathcal{K}_{0,2}=\mathcal{H}_{0,2}-\frac{\mu}{2a}\frac{R_\Earth^4}{p^4}\frac{3\eta}{32}\delta^2\left[(3e^2+2)(35s^4-40s^2+8)-10(7s^2-6)s^2e^2\cos2\omega \right].
\end{equation}
Because both models match to the first order of $J_2$, they also share the transformation that removes the short-period terms up to this order. Its computation is standard and hence the details on the derivation are omitted. For completeness, the corresponding corrections are summarized in Appendix \ref{a:shortperiod}.
\par

The long-term dynamics is qualitatively the same for both the $J_2$ and $J_4$ model, and Vinti's truncated model. That is, on average, all the elliptic orbits have rotating perigees for low inclinations; the mean circular orbits undergo a double pitchfork bifurcation close to the critical inclination, after which elliptic orbits with either rotating or oscillating perigees exist in different regions of the reduced phase space. The interested reader is referred to the illustrative Fig.~5 of \citet{CoffeyDepritDeprit1994} for details on the $J_2$ and $J_4$ model, and to Fig.~2 of \citet{Lara1996} for the truncation $\mathcal{O}(J_2^2)$ of Vinti's dynamical model.
\par

\begin{figure}[htbp] \center
\includegraphics[scale=0.63]{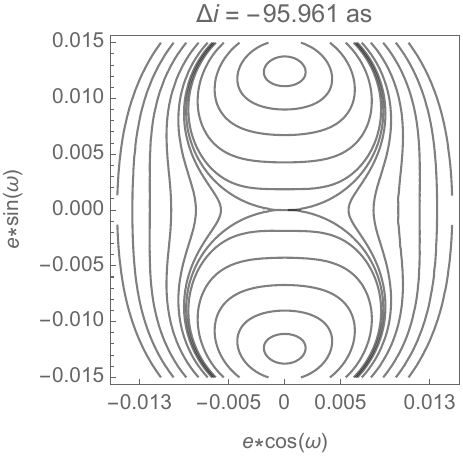} \hfil \includegraphics[scale=0.63]{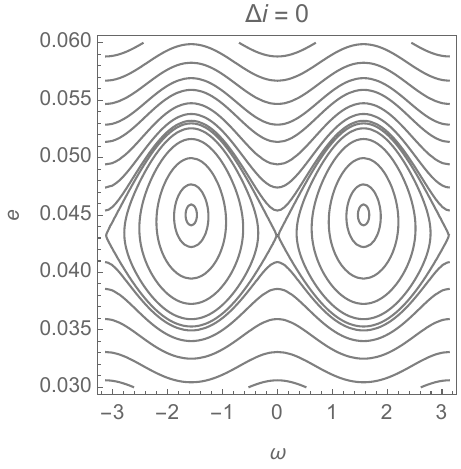}
 \\[0.5ex]
\includegraphics[scale=0.63]{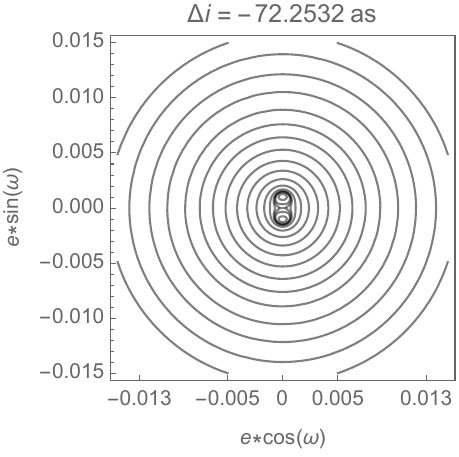} \hfil \includegraphics[scale=0.63]{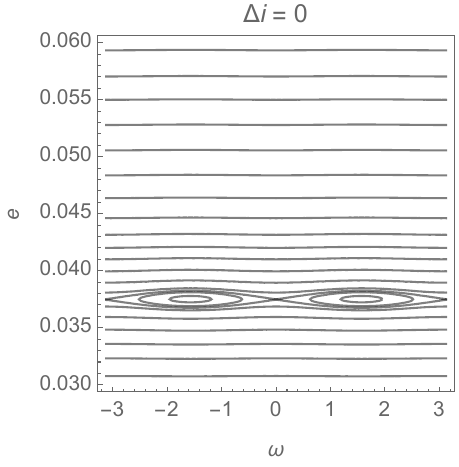}
\caption{Mean phase space diagrams. Above: $J_2$ and $J_4$ model. Below: $\mathcal{O}(J_2^2)$ truncation of Vinti's dynamical model}
\label{f:bifIvsVi}
\end{figure}

On the other hand, the naturalization process that turned $\mathcal{K}$ formally into $\mathcal{H}$ has introduced significant quantitative differences in the mean dynamics, shrinking the regions of orbits with librating mean periapsis by about one order of magnitude. This is illustrated in Fig.~\ref{f:bifIvsVi}, where the reduced $(g',G')$ dynamics is described in terms of the eccentricity and the argument of the perigee by means of polar (left column) and cylindrical maps (right column). The plots are depicted for the value $J_2(R_\Earth/a)^2=10^{-3}$ in order to enhance the effect of the perturbation, and the label $\Delta{i}$ denotes the difference between the critical inclination and the mean, constant inclination of the circular orbits. 
\par

The eccentricity vector diagrams in the left column of Fig.~\ref{f:bifIvsVi} correspond to the approximate moment in which the second pitchfork bifurcation takes place. The bifurcation occurs for different inclinations of the mean circular orbits depending on the model, although it is always close to the critical inclination. Hence, these graphics illustrate the different sizes of the regions of orbits with oscillating perigees of the $J_2$ and $J_4$ model, and Vinti's truncated model (upper and lower plots respectively). As clearly noticed, the former is roughly one order of magnitude larger than the later.
\par

The cylindrical maps in the right column of Fig.~\ref{f:bifIvsVi} correspond to the case in which the circular orbits of each different model are critically inclined. 
They confirm that the stripes in which the orbits with mean oscillating live remain of different orders of magnitude. They also illustrate that the mean dynamics outside that stripe remains much more regular, almost flat, in the case of Vinti's model (right-lower plot). While the consequent removal of the long-period terms will exclude the regions of orbits with librating perigees from the secular solution, the higher regularity of the mean dynamics resulting from Vinti's model should provide good expectations of having a secular solution that more faithfully represents the regions close to the strip.
\par

\subsection{Secular terms}

The first order terms of both $\mathcal{H}^\ddag$ and $\mathcal{K}'$ are independent of the argument of the perigee, and hence they remain unaltered by the second Lie transformation. As a consequence, the first order term of the generating function of the new transformation cannot be determined until the next order of the perturbation approach.
\par

At the second order, the simple averaging over the argument of the perigee yields the new Hamiltonian term
\begin{align} \label{Hp02}
\mathcal{H}^\ddag_{0,2}= &-\frac{\mu}{2a}\frac{R_\Earth^4}{p^4}\eta\frac{3}{16}\left[
5(21s^4-28s^2+8)+2(3s^2-2)^2\eta-2(25s^4-32s^2+8)\eta^2 \right],
\end{align}
for Vinti's truncated Hamiltonian, and
\begin{align} \label{Kp02}
\mathcal{K}'_{0,2}= &\;\mathcal{H}'_{0,2}-\frac{\mu}{2a}\delta^2\frac{R_\Earth^4}{p^4}\eta\frac{3}{32}(35s^4-40s^2+8)(3e^2+2)
\end{align}
in the case of the $J_2$ and $J_4$ Hamiltonian.
\par

The first-order term of the new generating function, say $\mathcal{U}_1$, can now be computed;  
refer to \citep{Brouwer1959} or section 6.3 of \citep{Lara2021} for additional details. For the $J_2$ and $J_4$ problem, it is solved from
\begin{equation} \label{homological}
\{\mathcal{U}_1,\mathcal{K}_{0,1}\}=\frac{1}{2}(\mathcal{K}_{0,2}-\mathcal{K}'_{0,2}),
\end{equation}
in which $\mathcal{K}_{0,1}=\mathcal{H}_{0,1}$ is replaced from Eq.~(\ref{H01}) to yield
\begin{equation} \label{longhomological}
\{\mathcal{U}_1,\mathcal{K}_{1,0}\}=\frac{\partial\mathcal{U}_1}{\partial{g}'}\frac{\partial\mathcal{K}_{1,0}}{\partial{G'}}\equiv-n_\omega\frac{\partial\mathcal{U}_1}{\partial{g}'},
\end{equation}
where
\begin{equation}
n_\omega 
=n\frac{R_\Earth^2}{p^2}\frac{3}{4}(5s^2-4),
\end{equation}
which vanishes at the critical inclination.
\par

Then, plugging Eq.~(\ref{longhomological}) into Eq.~(\ref{homological}), in which $\mathcal{K}_{0,2}$ and $\mathcal{K}'_{0,2}$ are replaced from Eqs.~(\ref{Kprime}) and (\ref{K02}), respectively, Eq.~(\ref{homological}) is solved by indefinite integration, to obtain
\begin{equation} \label{U1J2J4}
\mathcal{U}_1 
=-G'\frac{R_\Earth^2}{p^2}\frac{1}{32}\left(7\delta^2-\frac{2\delta^2}{5s^2-4}+4\right)s^2e^2\sin2\omega,
\end{equation}
and hence the critical divisor.
\par

Regarding Vinti's truncated model, the generating function is analogously computed replacing $\mathcal{K}_{0,1}$, $\mathcal{K}_{0,2}$, and $\mathcal{K}'_{0,2}$ in Eq.~(\ref{homological}), by $\mathcal{H}_{0,1}$, $\mathcal{H}_{0,2}$, and $\mathcal{H}^\ddag_{0,2}$, respectively, as given by Eqs.~(\ref{H01}), (\ref{H02}) and (\ref{Hp02}). Now, 
\[
\mathcal{H}_{0,2}-\mathcal{H}^\ddag_{0,2}=-\frac{\mu}{2a}\frac{R_\Earth^4}{p^4}\eta\frac{3}{4}(5s^2-4)s^2e^2\sin2\omega,
\]
and the factor $5s^2-4$ cancels the critical divisor resulting from the indefinite integration, yielding
\begin{equation} \label{U1vinti}
\mathcal{U}_1=-G'\frac{R_\Earth^2}{p^2}\frac{1}{8}s^2e^2\sin2\omega,
\end{equation}
thus completely avoiding the trouble of small divisors.
\par

\section{A more significant surrogate} \label{s:morsig} 

While the previous solution achieves the goal of avoiding the trouble of small divisors, its secular frequencies depart from the correct values in that they no longer depend on $J_4$. The origin of the discrepancy is easily identified with the secular terms introduced in the naturalization step, as follows from visual inspection of $\mathcal{S}$ in Eq.~(\ref{secular}), and the last term of $\mathcal{K}'_{0,2}$ in Eq.~(\ref{Kp02}). Note that the contribution of the latter to the secular dynamics must be multiplied by the coefficient ${1}/{2!}$ of the perturbation algorithm. Therefore, computing an infinitesimal contact transformation from a generating function that ignores $\mathcal{S}$ in Eq.~(\ref{Wnatu}) will keep the proper secular terms of the $J_2$ and $J_4$ model. 
\par

That is, the naturalization step is carried out replacing Eq.~(\ref{Wnatu}) by the new generating function
\begin{equation} \label{Wmixed}
\mathcal{W}=-\Theta{J}_2^2\delta^2\frac{R_\Earth^4}{p^4}\left(\mathcal{M}\theta+\mathcal{P}\right),
\end{equation}
in which $\mathcal{M}$ and $\mathcal{P}$ are still given by Eqs.~(\ref{natumixed}) and (\ref{natuperi}), respectively. The infinitesimal contact transformation is yet contaminated with terms that may grow unbounded, but the fact that it does no longer depends on purely secular terms makes it more appealing than in the previous case. In particular, the new infinitesimal contact transformation takes the form
\begin{equation}
\{\xi,\mathcal{W}\}=J_2^2\delta^2\frac{R_\Earth^4}{p^4}\Delta_\xi.
\end{equation}
For the polar-nodal variables, straightforward computations yield
\begin{align} \nonumber
\Delta_r=&-p\theta\frac{15}{32}(7s^2-6)s^2(\sigma\cos2\theta -\kappa\sin2\theta) 
-p\frac{1}{256}\left\{6(2\kappa^2+4\sigma^2+3\kappa+6)(35s^4 \right. \\ \nonumber
& -40s^2+8)
-2(26\kappa^2-12\sigma^2+35\kappa+14)s^4\cos4\theta +(64\kappa+35)\sigma{s}^4\sin4\theta
\\ \label{dr}
& \left. -2 \left(7 s^2-6\right) s^2 \left[ (8\kappa^2-48\sigma^2+15\kappa-40)\cos2\theta +(64\kappa+45)\sigma\sin2\theta\right]   \right\},
\\ \nonumber 
\Delta_\theta=&-\theta\frac{15}{64}\big\{4\left[2(14s^4-16s^2+3)\kappa-s^2(7s^2-6)\right]\sigma\sin2\theta
+\big[(49s^4-58s^2+12)\kappa^2 \\ \nonumber
&-(63s^4-70s^2+12)\sigma^2-4(7s^2-6)s^2\kappa\big]\cos2\theta \big\}
+\frac{1}{512}\big\{8\big[18(35s^4-50s^2 \\ \nonumber
&+16)\kappa^2+8(70s^4-95s^2+28)\sigma^2+6(140s^4-205s^2+68)\kappa+135(21s^4 \\ \nonumber
&-28s^2+8)\big]\sigma
+8\big[32(14s^4-16s^2+3)\sigma^2-8(21s^4-26s^2+6)\kappa^2-2(154s^4 \\ \nonumber
&-192s^2+45)\kappa+5(301s^4-322s^2+48)\big]\sigma\cos2\theta -2\big[32(35s^4-46s^2+12)\kappa^3 \\ \nonumber
& +64(49s^4-58s^2+12)\kappa\sigma^2 +3(777s^4-1066s^2+300)\kappa^2 +20(399s^4 \\ \nonumber
&-470s^2+96)\kappa +(1939s^4-2382s^2+540)\sigma^2 +40(21s^4-34s^2+12)\big]\sin2\theta \\ \nonumber
&+8 s^2\big[8(2s^2-1)\sigma^2-26(3s^2-2)\kappa^2-2(44s^2-35)\kappa-35s^2+28\big]\sigma\cos4\theta \\ \nonumber
&+s^2\big[96(5s^2-4)\kappa^3-64(7s^2-4)\kappa\sigma^2 +(1139s^2-980)\kappa^2-(187s^2-140)\sigma^2 \\ \label{dz}
&+28(37s^2-32)\kappa+14(23s^2-20)\big]\sin4\theta \Big\},
\\ \nonumber
\Delta_\nu=&
-\theta\frac{15}{16}c(7s^2-3)\big[\left(\kappa^2-\sigma^2\right)\cos2\theta +2\kappa\sigma\sin2\theta\big]
+\frac{1}{128}c\big\{ 20(7s^2-4)\sigma(6\kappa^2 \\ \nonumber
&+4\sigma^2+9\kappa+18) 
 +2(7s^2-3)\big[2\sigma(16\sigma^2-8\kappa^2-15\kappa+40)\cos2\theta
 -(32\kappa^3 \\ \nonumber
&+64\kappa\sigma^2+75\kappa^2+45\sigma^2+160\kappa+40)\sin2\theta \big]
+4s^2(4\sigma^2-26\kappa^2-35\kappa-14) \\ \label{dnu}
&\times \sigma\cos4\theta 
+s^2(96\kappa^3-64\kappa\sigma^2+245\kappa^2-35\sigma^2+224\kappa+70)\sin4\theta \big\},
\\ \nonumber 
\Delta_R=&\; \theta\frac{\Theta}{r}\frac{15}{32}\frac{p}{r}s^2(7s^2-6)(\kappa\cos2\theta+\sigma\sin2\theta) 
-\frac{\Theta}{r}\frac{1}{256}\frac{p}{r}\big\{6(35s^4-40s^2+8)\sigma(4\kappa \\ \nonumber
&+3)
-2(7s^2-6)s^2\left[(48\kappa^2+32\sigma^2+75\kappa+80)\sin2\theta +(16\kappa+15)\sigma\cos2\theta\right] \\ \label{dvr}
&+(144\kappa^2-32\sigma^2+245\kappa+112)s^4\sin4\theta - 2(52\kappa+35)s^4\sigma\cos4\theta \big\},
\\ \nonumber 
\Delta_\Theta=&\; \theta\Theta\frac{15}{32}s^2(7s^2-6)\left[(\kappa^2-\sigma^2)\sin2\theta -2\kappa\sigma\cos2\theta\right]
+\Theta\frac{1}{128}\big\{ 4\sigma(26\kappa^2-4\sigma^2 \\ \nonumber
&+35\kappa+14)s^4\sin4\theta +(96\kappa^3-64\kappa\sigma^2+245\kappa^2-35\sigma^2+224\kappa+70)s^4\cos4\theta \\ \nonumber
& -(7s^2-6) s^2\big[ 2\sigma(16\sigma^2-8\kappa^2+15\kappa+40)\sin2\theta
+(32\kappa^3+64\kappa\sigma^2+105\kappa^2 \\ \label{dG}
& +15\sigma^2+160\kappa+40)\cos2\theta \big] \big\}.
\end{align}
\par

In the new variables, the new perturbation Hamiltonian only differs from Eq.~(\ref{H20}) in the component that gives rise to the secular term of concern. Namely,
\begin{align} \label{H20ast}
\mathcal{H}^*_{2,0}= &\;\mathcal{H}_{2,0}-\frac{\mu}{2a}\delta^2\frac{R_\Earth^4}{p^4}\left(\frac{a^2\eta}{r^2}\right)\frac{3}{32}\eta(35s^4-40s^2+8)(3e^2+2).
\end{align}
Obviously, since the new term is of the second order of $J_2$, it does not affect the short-period corrections, which remain the same as before.
\par

The corresponding Hamiltonian in mean elements only differs from Eq.~(\ref{Vinti22}) in the secular term. Namely,
\begin{align} \label{Hp02ast}
\mathcal{H}^*_{0,2}= &\;\mathcal{H}_{0,2}-\frac{\mu}{2a}\delta^2\frac{R_\Earth^4}{p^4}\frac{3}{32}\eta(35s^4-40s^2+8)(3e^2+2)
\end{align}
where $\mathcal{H}_{0,2}$ is the same previously given in Eq.~(\ref{H02}). The following removal of long-period terms preserves the last term on the right side of Eq.~(\ref{Hp02ast}), thus achieving the desired, same result as in Eq.~(\ref{Kp02}). 
\par

Because secular terms are not affected by the long-period elimination, the new generating function remains the same as Eq.~(\ref{U1vinti}), from which, to the first order, the long-period corrections are obtained from the straightforward evaluation of Poisson brackets, say
\begin{equation}
\{\xi,\mathcal{U}_1\}=J_2\Delta'_\xi.
\end{equation}
\par

On the other hand, the new theory is easily suplemented with the contribution of the third zonal harmonic since it does not involve critical divisors in the removal of long-period terms, and, to the precision of such kind of theory the short-period corrections are unaffected by second-order terms, cf.~\citep{Brouwer1959}. That is, the term $\mathcal{K}_{2,0}$ in Eq.~(\ref{geo2even}) is now supplemented with the additional terms
\begin{equation} \label{J3oscu}
\Delta_{J_3}=\frac{\mu}{r}\tilde{J}_3\frac{R_\Earth^3}{r^3}\frac{1}{2}\left[3(5s^2-4)s\sin\theta-5s^3\sin3\theta\right],
\end{equation}
with $\tilde{J}_3\equiv{J}_3/J_2^2$, which follow from the standard expansion of the zonal potential and the consequent scaling introduced by the perturbation arrangement. Thus, Eq.~(\ref{U1vinti}) turns into
\begin{equation}
\mathcal{U}_1=-G\frac{R_\Earth^2}{p^2}\frac{1}{8}
s^2e^2\sin2\omega
+G\tilde{J}_3\frac{R_\Earth}{p}\frac{1}{2}se\cos\omega.
\end{equation}
\par

In particular, for the usual set of semi-equinoctial, nonsingular variables for circular orbits, namely
\begin{equation}
\lambda=\ell+\omega, \qquad 
k=e\cos\omega, \qquad
q=e\sin\omega, 
\end{equation}
which are customarily complemented with $a$, $i$, and $\Omega$, it is readily obtained
\begin{align} \label{lalp}
\Delta'_\lambda=& \; \frac{1}{4}\frac{R_\Earth^2}{p^2}{k}q\Big[s^2\left(5+\frac{2\eta^2}{1+\eta}\right)-2\Big]
-\frac{1}{2}\tilde{J}_3\frac{R_\Earth}{p}k\left[s\left(\frac{\eta^2}{1+\eta}+2\right)-\frac{1}{s}\right] ,
\\ \label{klp}
\Delta'_k=&-\frac{1}{4}\frac{R_\Earth^2}{p^2}{k}\big[(4s^2-2)q^2-s^2(k^2-1)\big] 
+\frac{1}{2}\tilde{J}_3\frac{R_\Earth}{p}\frac{1}{s}kq(2s^2-1),
\\ \label{qlp}
\Delta'_q=&\; \frac{1}{4}\frac{R_\Earth^2}{p^2}{q}\big[(4s^2-2)k^2-s^2(q^2-1)\big]
-\frac{1}{2}\tilde{J}_3\frac{R_\Earth}{p}\frac{1}{s}\left[(s^2-1)k^2-s^2(q^2-1)\right],
\\
\Delta'_a=& \; 0, \\ \label{ilp}
\Delta'_i=&\; \frac{1}{4}\frac{R_\Earth^2}{p^2}(k^2-q^2)cs
+\frac{1}{2}\tilde{J}_3\frac{R_\Earth}{p}cq,
\\ \label{hlp}
\Delta'_\Omega=&\; \frac{1}{2}\frac{R_\Earth^2}{p^2}kqc
-\frac{1}{2}\tilde{J}_3\frac{R_\Earth}{p}\frac{c}{s}k.
\end{align}
Division by the sine of the inclination related to the third zonal harmonic would be removed, if desired, by choosing a different set of non-singular variables as, for instance, the Poincar\'e variables used by \citet{Lyddane1963}.
\par

The secular terms of the new $J_2$--$J_4$ solution are, therefore, the same of Brouwer's theory when ignoring the contribution of the fifth zonal harmonic from the latter. Regarding them, it is worth recalling that the secular frequencies are constant, and hence, the secular argument of the perigee either grows continuously, at positive or negative rate, or remains constant (frozen). Therefore, like in the case of intermediary solutions, the secular terms predict that, for given $a$ and $e$, all orbits are frozen at a certain inclination---yet now the inclination departs slightly from the critical value. They are suitably reformulated in the same set of variables, and are provided in Appendix \ref{a:secular} for the convenience of prospective readers.
\par


The previous discussion refers to the analytical solution, which, as it was already emphasized, does not cover the region of orbits with librating perigees. Regarding the mean dynamics, the changes introduced by the naturalization process make that the effects of the third zonal harmonic may now dominate over those of the fourth in some regions of phase space, thus clearly modifying the long-term behavior of the orbits with oscillating perigees.
\par

A complete bifurcation analysis of the mean dynamics falls out the scope of the current investigation. We are satisfied with illustrating the differences introduced by the naturalization with the help of some diagrams describing the mean dynamics of the eccentricity and the argument of the perigee.
\par

\begin{figure}
\center
\includegraphics[scale=0.65]{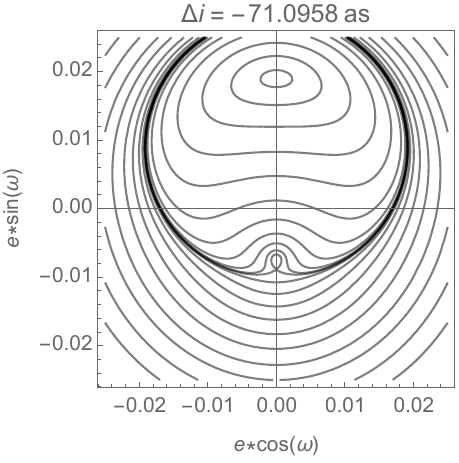} \hfil
\includegraphics[scale=0.65]{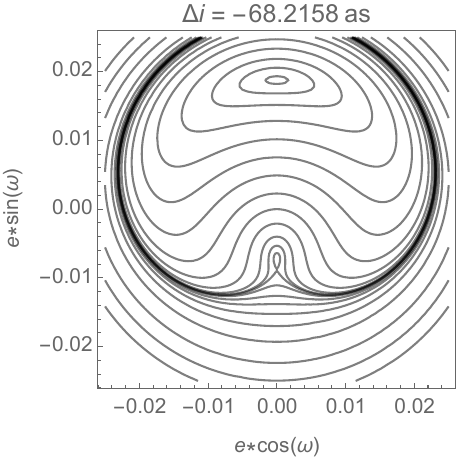}
\caption{$J_2$--$J_4$ model ($a=1.16R_\Earth$). Left: classical approach. Right: new theory. }
\label{f:bifJ3}
\end{figure}

As it is well known, the third zonal harmonic coefficient breaks the equatorial symmetry of the $J_2$ and $J_4$ problem; cf.~Fig.~(5) of \citet{CoffeyDepritDeprit1994}. For orbits below the critical inclination, the new approach only produces small quantitative changes in the long term behavior. This is illustrated in Fig.~\ref{f:bifJ3} for the case in which both theories, the classical (left plot) and the new one (right plot), have just undergone the first bifurcation of frozen orbits, which happens for slightly different values of $i_\mathrm{circular}$. The bifurcation emerges in the form of a teardrop, with one stable and one unstable frozen orbits, both with the mean perigee fixed at $270$ degree, the stable frozen orbit having the smaller eccentricity. The regions of orbits with librating and rotating perigees remain of similar sizes in both cases, with a slightly larger amplitude of the perigee oscillations in the case of the new theory.
\par

\begin{figure}
\center
\includegraphics[scale=0.63]{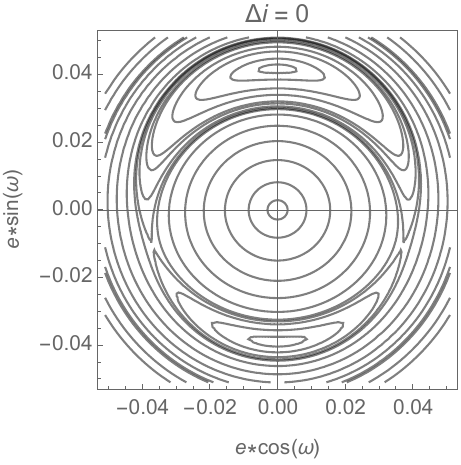} \hfil \includegraphics[scale=0.63]{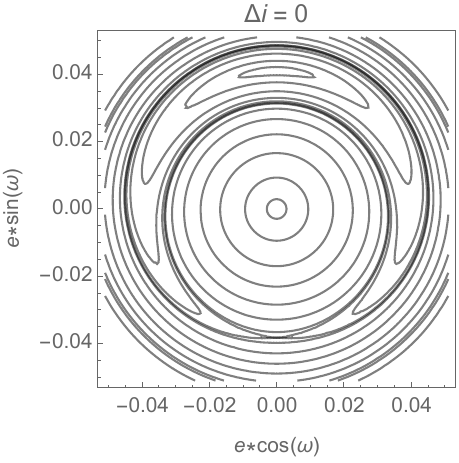} \\
\includegraphics[scale=0.63]{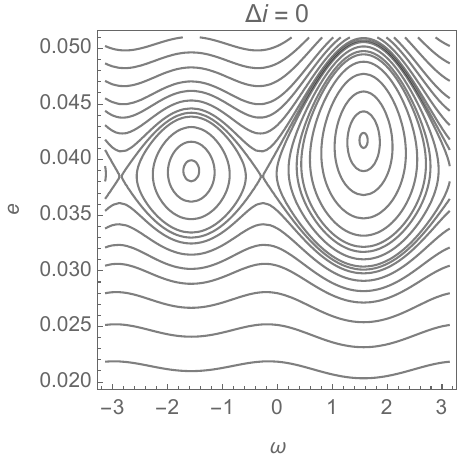} \hfil \includegraphics[scale=0.63]{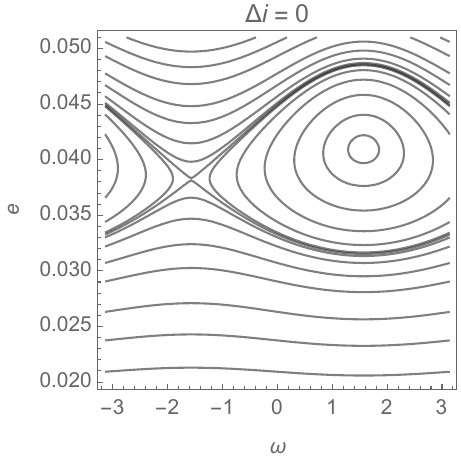}
\caption{$J_2$--$J_4$ model ($a=1.16R_\Earth$). Left: classical approach. Right: new theory. }
\label{f:2bif}
\end{figure}

Things change radically, nevertheless, for slightly higher inclinations of the circular orbits, as shown in Fig.~\ref{f:2bif} for the case in which $i_\mathrm{circular}$ matches the critical inclination value. The bifurcated stable frozen orbits reduced their mean eccentricity in both theories, whereas the mean eccentricity of the other bifurcated orbits increased their value. However, within the classical theory, the unstable frozen orbit changed to stability in a pitchfork bifurcation, in which two new unstable frozen orbits emerged, whose perigees depart symmetrically from the $270$ degree direction (left plot of Fig.~\ref{f:2bif}). On the contrary, there are not qualitative changes in the mean dynamics undergone by the new theory, which maintains essentially the $J_2$--$J_3$ dynamics. The plots in the second row of Fig.~\ref{f:2bif} illustrate the same situation by means of cylindrical $\omega$-$e$ maps, focusing on the different regions of orbits with oscillating perigees. While the region of libration are of similar amplitude in both cases, the orbits with mean rotating perigee evolve more regularly in the case of the new theory (right-lower plot).
\par

\begin{figure}
\center
\includegraphics[scale=0.68]{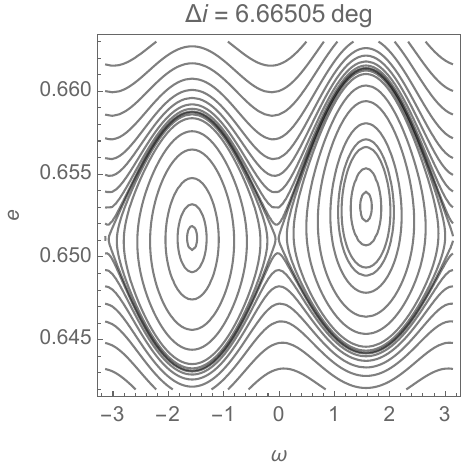} \hfil \includegraphics[scale=0.68]{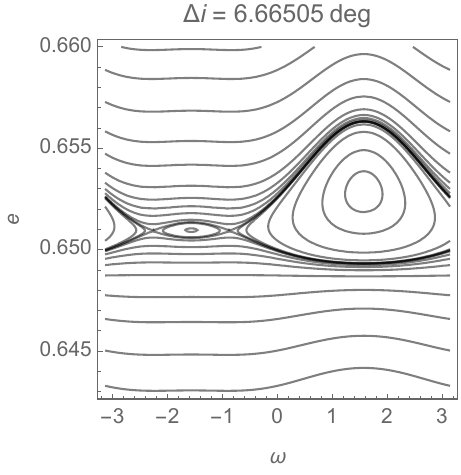}
\caption{$J_2$--$J_4$ model ($a=1.16R_\Earth$). Left: classical approach. Right: new theory. }
\label{f:2bifIntJ3}
\end{figure}

The second bifurcation also occurs in the case of the new theory, but it is separated by several degrees from the critical inclination. This is shown in Fig.~\ref{f:2bifIntJ3}, where $i_\mathrm{circular}\approx73$ degree. The right plot shows that the bifurcation just happened for the new theory. In spite of that, the behavior is more regular with the new theory, and the size of the regions of elliptic orbits with oscillating perigee, remains notably smaller than in the case of the classical theory shown in the left plot of the figure. It must be noted that, due to the high eccentricity and low altitude of the orbits illustrated in Fig.~\ref{f:2bifIntJ3}, they would impact the surface of the Earth. However, the case presented is not intended for practical application, but to illustrate the differences on the mean dynamics resulting for each different approach.
\par

%


\section{Example propagations}


The performance of the new theory is illustrated with two different examples of orbits close to the critical inclination. Namely, a low- and a high-eccentricity orbit. The accuracy of the analytical solution is assessed  by comparison with reference orbits that are integrated numerically in Cartesian coordinates. Checking that the numerical solution preserves 13 digits at least in each integral of the problem---the energy and the third component of the angular momentum vector---is considered a satisfactory quality test of the true orbit. 

For increased accuracy of the analytical propagations, the initialization of the perturbation solution from given initial conditions is improved with the energy calibration of the semimajor axis, in this way notably reducing the abnormal growth of along-track errors \citep{BreakwellVagners1970,Lara2020arxiv}.
\par

The values $\mu=398600.4415\,\mathrm{km^3/s^2}$, $R_\Earth=6378.1363$ km, $J_2=1.08264\times10^{-3}$, $J_3=-2.53247\times10^{-6}$, and $J_4=-1.61997\times10^{-6}$, are assigned in the simulations to the parameters defining the potential model---in rough agreement with the unnormalized values of the GRACE gravity model GGM02 \citep{Tapleyetal2005}.

The first test case is for a critically inclined, Molnya-type orbit, which is propagated analytically for an interval of 4 weeks. The initial osculating elements are $a=26554$ km, $e=0.72$, $i=63.425^\circ$, $\Omega=5.73^\circ$, $\omega=280^\circ$, and $M=0$. From these values, the secular terms are initialized following the sequence in Table~\ref{t:Molnya}, where ``naturalized'' refers to the values obtained after the infinitesimal contact transformation defined in Eqs.~(\ref{dr})--(\ref{dG}). The 10-digit accuracy provided in the tables should be enough for an interested reader to reproduce the obtained results. Recall that $N$ is an integral of the original dynamics.
\par

\begin{table}[htbp]  
\caption{ Molnya-type orbit. Initialization sequence for the analytical theory. Distances are km, velocities km/s, momenta km${}^2$/s, and angles radians.
\label{t:Molnya}
} 
\begin{tabular}{@{}rlrrrrc@{}} \noalign{\smallskip}
\multicolumn{2}{c}{Osculating} & \multicolumn{1}{c}{Naturalized} & \multicolumn{1}{c}{Mean} & \multicolumn{2}{c}{Secular} & \multicolumn{1}{c}{Calibrated}\\
\hline\noalign{\smallskip}
$r=$      & $7435.12$      & $ 7435.119763$ & $ 7434.999020$ & $\lambda=$ & $  4.887011931$ & \\
$\theta=$ & $4.886921906$  & $ 4.886921948$ & $ 4.886930581$ & $k=$       & $ 0.1252997492$ & \\
$\nu=$    & $0.1000073661$ & $0.1000072391$ & $0.1000678110$ & $q=$       & $-0.7103138281$ & \\
$10^8R=$  & $0.0$          & $-2.289903912$ & $-60842.44996$ & $a=$       & $  26653.67292$ & $26653.76699$ \\
$\Theta=$ & $71396.52959$  & $ 71396.52738$ & $ 71417.77113$ & $i=$       & $  1.106951853$ & \\
$N=$      & $31940.58668$  &                &                & $\Omega=$  & $ 0.1000366312$ & \\
\noalign{\smallskip}\hline
\end{tabular}
\end{table}

The formulas for computing the mean to osculating transformation only involve the effects of the second zonal harmonic, and are standard. The use of polar variables is encouraged for a compact representation of the short-period corrections \citep{Kozai1962,Izsak1963}, which we used directly in the practicable form given in \citep{HautesserresLaraCeMDA2017}. For completeness, they are summarized in Appendix \ref{a:shortperiod}.
\par

\begin{figure}[htbp]
\includegraphics[scale=0.66]{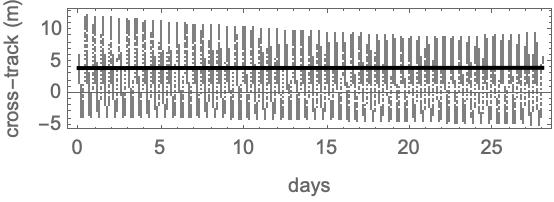} \;
\includegraphics[scale=0.66]{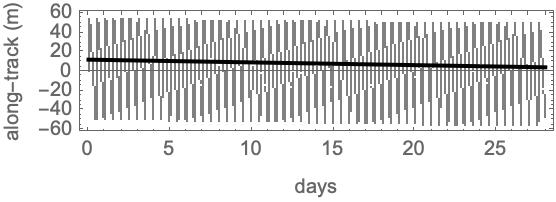} \\
\includegraphics[scale=0.66]{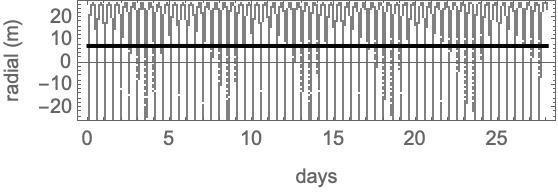} \;
\includegraphics[scale=0.66]{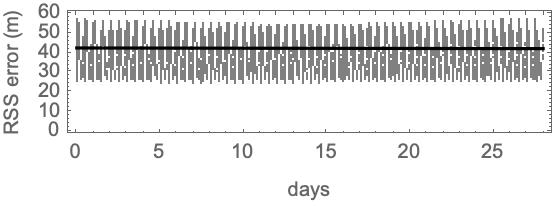}
\caption{Molnya-type example orbit. Position errors of the new analytical solution.}
\label{f:intrinsicMolnya4w}
\end{figure}

The intrinsic errors obtained with the analytical solution are displayed in Fig.~\ref{f:intrinsicMolnya4w}. There is no appreciably secular trend in the errors obtained with the new theory, the largest one of about 30 cm/day corresponding to the case of along-track errors  (upper-right plot of Fig.~\ref{f:intrinsicMolnya4w}). That the new theory closely captures the dynamics is clearly observed in the proper modeling of the evolution of the orbital plane, as noted from the time-history of the cross-track errors (upper-left plot of Fig.~\ref{f:intrinsicMolnya4w}), which are only affected of a short-period modulation of about 15 m of amplitude. But this is for the insufficient modulation of the short-period corrections due to their truncation to the first order of $J_2$, which is an unavoidable limitation of both analytical solutions. When the position errors are combined into a single measure, we obtain the root sum squared (RSS) errors displayed in bottom-right plot of Fig.~\ref{f:intrinsicMolnya4w}. The time history of the RSS errors shows a negligible secular growth of about one cm by day, which is modulated with short-period oscillations of just a few tens of meter of amplitude.
\par

\begin{figure}[htbp]
\includegraphics[scale=0.66]{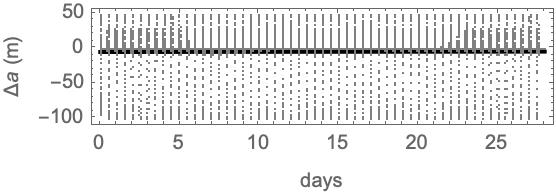} \;
\includegraphics[scale=0.66]{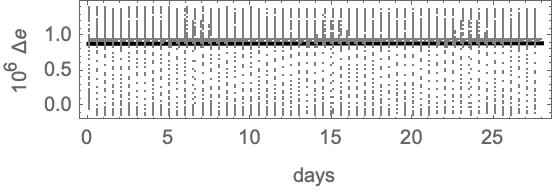} \\
\includegraphics[scale=0.66]{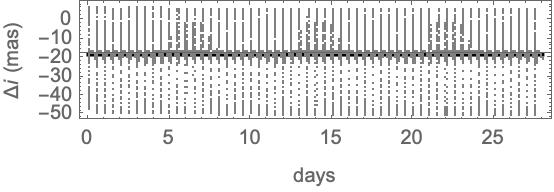} \;
\includegraphics[scale=0.66]{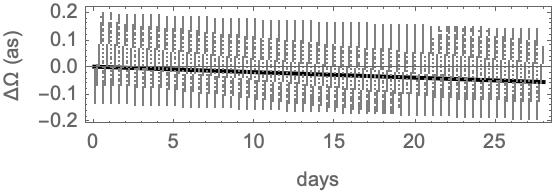} \\
\includegraphics[scale=0.66]{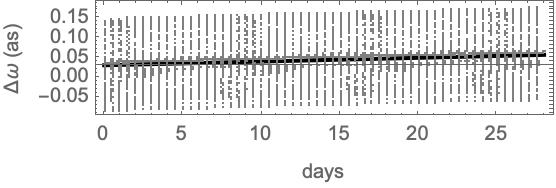} \;
\includegraphics[scale=0.66]{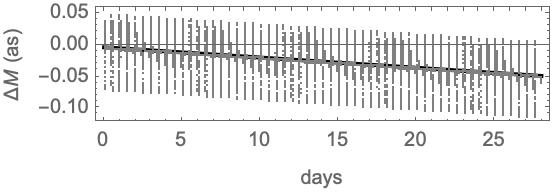} 
\caption{Critically inclined, Molnya-type example orbit. Time history of the errors of the orbital elements.}
\label{f:DifMolnya4w}
\end{figure}

The merits of the new solution are further illustrated in Fig.~\ref{f:DifMolnya4w}, which displays the time histories of the errors of the orbital elements. The errors of the semimajor axis mostly average to zero, with a negligible secular growth of a bout 2 cm/day, which was expected as long as the semimajor axis is not affected of terms that may grow unbounded by the naturalization. The errors of the eccentricity and inclination have negligible secular trends, and those of the other angles are certainly small, of the order of the milli arc second by day.
\par

While the mixed secular-periodic terms do not produce harmful, observable effects in this interval of about 60 orbital periods, for sure they will disclose in much longer intervals.
\par

The second test is for a much lower, low-eccentricity orbit, whose inclination departs somewhat from the critical value. It is now analytically propagated for about one week and a half. This time interval roughly corresponds now to 130 orbital periods, a span that should allow for a more demanding test on the effects of the unbounded terms introduced by the naturalization. The initial, osculating elements are $a=7707.27$ km, $e=0.09$, $i=62.5^\circ$, $\Omega=90^\circ$, $\omega=180^\circ$, and $M=0$. The sequence for the initialization of the secular terms is presented in Table~\ref{t:Topex}.
\par

\begin{table}[htb]  \tabcolsep 6pt  \centering 
\caption{ Low-eccentricity orbit. Initialization sequence for the analytical theory. Distances are km, velocities km/s, momenta km${}^2$/s, and angles radians.
\label{t:Topex}
}
\begin{tabular}{@{}rlrrrrc@{}} \noalign{\smallskip} 
\multicolumn{2}{c}{Osculating} & \multicolumn{1}{c}{Naturalized} & \multicolumn{1}{c}{Mean} & \multicolumn{2}{c}{Secular} & \multicolumn{1}{c}{Calibrated}\\
\hline\noalign{\smallskip}
$r=$      & $7013.615700$ & $7013.614926$  & $7010.990861$  & $\lambda=$ & $ 3.141687487$ & \\
$\theta=$ & $3.141592654$ & $3.141592671$  & $3.141592672$  & $10^2k=$   & $-8.933507464$ & \\
$\nu=$    & $1.570796327$ & $1.570796333$  & $1.570796333$  & $10^4q=$   & $-8.643997830$ & \\
$10^8R=$  & $0$           & $9.560507141$  & $9.570703947$  & $a=$       & $ 7698.648948$ & $7698.644923$ \\
$\Theta=$ & $55201.78478$ & $55201.77743$  & $55174.28440$  & $i=$       & $ 1.090570683$ & \\
$N=$      & $25489.34757$ &                &                & $\Omega=$  & $ 1.570841779$ & \\
\noalign{\smallskip}\hline
\end{tabular}
\end{table}

\begin{figure}[htbp]
\includegraphics[scale=0.66]{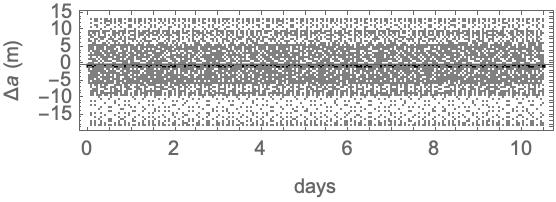} \;
\includegraphics[scale=0.66]{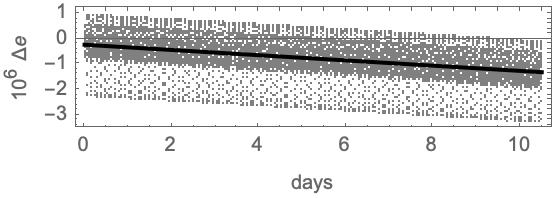} \\
\includegraphics[scale=0.66]{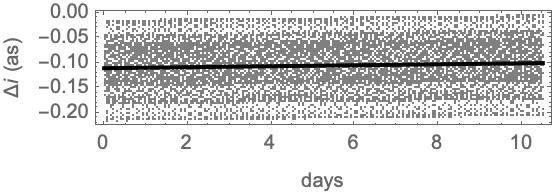} \;
\includegraphics[scale=0.66]{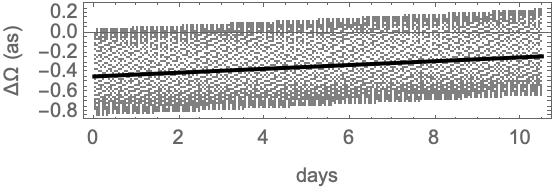} \\
\includegraphics[scale=0.66]{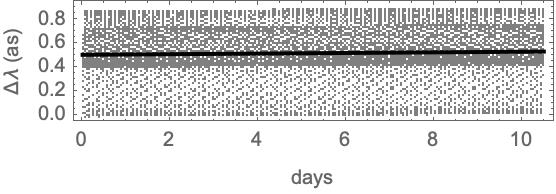} \;
\includegraphics[scale=0.66]{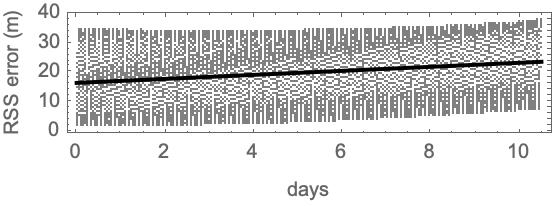} 
\caption{Low-eccentricity example orbit. Time history of the errors of the orbital elements obtained.}
\label{f:DifTopex10d}
\end{figure}

The time histories of the errors are now summarized in Fig.~\ref{f:DifTopex10d}. In every case, a linear fit (the black thick line) is superimposed to the time history to highlight the secular trend. Now, the errors of the eccentricity show a clear secular trend, as expected for the bigger effects of the perturbations due to the lower altitude of the orbit. The general behavior is otherwise analogous to the previous example. Now, due to the low eccentricity, the time history of the errors of the mean longitude $\lambda$ replace those of the mean anomaly and the argument of the perigee (bottom-left plot of Fig.~\ref{f:DifTopex10d}). As shown in the bottom-right plot of Fig.~\ref{f:DifTopex10d}, the position errors grow now at a much higher rate than in the previous example, slightly above half meter day. This is again due to the more perturbed dynamics, on the one hand, but also for the effects of the unbounded terms, on the other hand. The latter start to become apparent, say, after a week---or about 80 orbital periods---and would certainly make the predictions invalid eventually, in much longer time intervals.
\par

Needless to say, the new theory does not provide any advantage in the analytical computation of orbits whose inclinations depart enough from the critical value, a case in which the classical theory does not suffer any convergence problem. In practice, the new solution would replace the classical approach only when needed, typically for inclinations within 61 and 66 deg, cf.~\citep{CoffeyNealSegermanTravisano1995}, thus allowing for a much more accurate representation of the solution close to the critical inclination.
\par

\section{Possible extensions to higher orders} \label{s:higher}

Beyond first-order analytical solutions in Brouwer's style, it seems natural querying whether or not the trade off between small divisors and unbounded corrections can be extended to higher-order theories. The answer is positive, but with some precisions regarding the way of achieving it. More precisely, a naturalization based on a higher-order truncation of Vinti's dynamical system can certainly be applied. It would remove the critical divisors associated to the higher-order contribution of the even zonal harmonics, which appear both in the higher-order corrections and the higher-order secular terms. Needless to say, however, Vinti's dynamical system does not consider the odd zonal harmonics of the Geopotential.
\par

The fact that the third zonal harmonic is free from small divisors in Brouwer's first-order theory was crucial to applying the proposed naturalization. But this is not at all the case of higher-order theories, where the third zonal harmonic is affected of small divisors in the long-period corrections of the second order as well as in the third-order secular terms, cf.~\citep{CoffeyNealSegermanTravisano1995}. Additionally, the contribution of the fifth zonal harmonic, where critical divisors are concomitant to the long-period corrections as soon as at the first order, was ignored for the smallness of the Earth's $J_5$ coefficient. In these cases, it is remarked, the small divisors cannot be removed by a naturalization based on Vinti's dynamical system.
\par

The fact that elimination of divisors related to the critical inclination cannot be generally extended with the naturalization process presented here does not mean that the removal cannot be achieved. However, it would require the nontrivial task of identifying all the undesired terms of the zonal potential that give rise to the small divisor issue. A possible way of doing that would be approaching by Lie transforms the removal of only the ``nonessential'' long-period terms, cf.~\citep{Lara2021IAC}. Then, the remaining, essential long-period terms---which would introduce the small divisors in the classical approach---would be directly removed by a contact transformation. The latter, of course, should spring unbounded corrections. On the other hand, because the composition of Lie transformations is canonical \citep{Henrard1970Lie}, splitting \citep{Deprit1981}, composing \citep{Lara2020arxiv,LaraFantinoSusantoFlores2024}, or altering the order in which the different Lie transformations are applied \citep{AlfriendCoffey1984,CidFerrerSeinEchaluce1986,AbadSanJuanGavin2001} is mostly a matter the user's choice. Therefore, once the harming terms have been identified, the perturbation theory could be approached after a preliminary simplification that introduces the opposite of the undesired terms. In this way, the undesired inclination-resonant terms would cancel out in the following construction of the perturbation theory by Lie transforms or by any other perturbation approach. 
\par

Finally, it is worth mentioning that carrying out a preliminary simplification of the original Hamiltonian model based on the original ``elimination of the parallax'' \citep{Deprit1981}, or analogous simplifications inspired by it \citep{CoppolaPalacian1994,Lara2019CMDA}, could provide supplementary advantages regarding the implementation or evaluation of an operational code. 
\par

In summary, the design of a general method for dealing properly with all the critical divisors that may appear either in a higher-order perturbation theory, or when dealing with higher-fidelity models including the effects of additional zonal harmonics, as well as its possible effectivity for accurately propagating orbits with circulating periapsis close to the critical inclination, is a nontrivial problem that requires a careful study of the mentioned possible lines of attack. While it falls out of the scope of the current investigation, the topic is of undeniable interest, and could be approached as future research.
\par

\section{Conclusions}

That the naturalization of Vinti's dynamical system can be used to remove annoying divisors related to the critical inclination of $63.4^\circ$ from classical, first-order, analytical perturbation solutions, has been demonstrated. The technique is inspired by Deprit's original naturalization of common intermediaries of the main problem of the artificial satellite, and shares analogous merits and drawbacks. In particular, the proposed infinitesimal transformation bears secular, as well as mixed secular-periodic terms. Undeniably, this is a weak point of such naturalization. However, that the unbounded terms are confined to the second order of the Earth's oblateness coefficient is certainly good news, as it is also the fact that these kinds of terms do not corrupt the transformation of the semimajor axis. Moreover, beyond Vinti's dynamical system, the naturalization has been polished so that it only begets mixed secular-periodic terms, thus further delaying the appearance of harmful effects. In this sense, the new approach acts as a typical medical remedy. Namely, while the medicine is intended to fix a major health issue (small divisors), it may also cause adverse reactions (mixed terms) as a minor side effect of the healing.
\par

That existing analytical perturbation solutions based on total Hamiltonian reduction do not cover the case of orbits with oscillating perigees is well known. It is also the case of the new approach. However, because the new solution smooths the mean phase space close to the boundary separating the regions of rotating- and librating-perigee orbits, it avoids the usual lack of convergence of the series representing the long-period oscillations close to the critical inclination. Hence it better approaches the actual dynamics in that region. Needless to say, this is true only in the times in which the effects of the mixed terms remain buried by the errors stemming from the unavoidable truncation of perturbation solutions.
\par

That the proposed approach is constrained to the particular theories discussed (first order, and dealing only with the more significant harmonics of the Earth's gravitational potential) stems from the physical nature of Vinti's dynamical system. Possible avenues for extending the method to the case of zonal models of higher fidelity, or to progress to higher orders of the perturbation approach, have been outlined. However, devising a general procedure for identifying the inclination resonant terms of an arbitrary truncation of the zonal potential does not seem a trivial task. If considered of practical application, it would deserve further separate efforts.

\subsection*{Acknowledgements}
A preliminary version of this research was presented in the 
7th International Workshop on Key Topics in Orbit Propagation Applied to Space Situational Awareness (KePASSA2026), held in Toulon, France, May 27-29, 2026.


\begin{appendices}

\section{Secular terms} \label{a:secular}

The secular terms of the $J_2$--$J_4$ solution at a given time $t$ are computed from $a=a_0$, $i=i_0$, and
\begin{align} \label{secute}
\lambda=& \; \lambda_0+n_\lambda{t}, \\
k=& \; k_0\cos(n_\omega{t})-q_0\sin(n_\omega{t}), \\
q=& \; q_0\cos(n_\omega{t})+k_0\sin(n_\omega{t}), \\
\Omega=& \; \Omega_0+n_\Omega{t},
\end{align}
in which $\lambda_0$, $k_0$, $q_0$, $a_0$, $i_0$, and $\Omega_0$, denote the secular terms corresponding to the initial conditions of the analytical propagation. The involved secular frequencies are
\begin{align} \label{secula}
n_\lambda=& \; n-nJ_2\frac{3}{4}\frac{R_\Earth^2}{p^2}\left[\eta(3s^2-2)+5s^2-4\right]
+nJ_2^2\frac{3}{128}\frac{R_\Earth^4}{p^4}\sum_{j=0}^3\left(\Lambda_{1,j}-5\delta^2\Lambda_{2,j}\right)\eta^j, 
\\ \label{secuom}
n_\omega=& -nJ_2\frac{3}{4}\frac{R_\Earth^2}{p^2}(5s^2-4)
+nJ_2^2\frac{3}{128}\frac{R_\Earth^4}{p^4}\sum_{j=0}^2\left(\Gamma_{1,j}-5\delta^2\Gamma_{2,j}\right)\eta^j, 
\\ \nonumber
n_\Omega=& -nJ_2c\frac{3}{2}\frac{R_\Earth^2}{p^2}
+nJ_2^2c\frac{3}{32}\frac{R_\Earth^4}{p^4}\big[ 2(6\eta+35)(3s^2-2)-4\eta^2(25s^2-16) \\ \label{secuh}
&-5\delta^2(3\eta^2-5)(7s^2-4)\big],
\end{align}
where $n$ denotes the mean motion in secular variables, and the inclination polynomials $\Lambda_{i,j}$, and $\Gamma_{i,j}$, in Eqs.~(\ref{secula}) and (\ref{secuom}) are listed in Table~\ref{t:incpoly}.
\par

\begin{table}[htb]  \tabcolsep 3pt  \centering 
\caption{ Inclination polynomials of the secular terms.
\label{t:incpoly}
} 
\begin{tabular}{@{}lcccc@{}} \hline\noalign{\smallskip}
$j$ & $\Lambda_{1,j}$ & \multicolumn{1}{c}{$\Lambda_{2,j}$} & $\Gamma_{1,j}$ & $\Gamma_{2,j}$ \\
\hline\noalign{\smallskip}
$0$ & $\Gamma_{1,0}$  & $\Gamma_{2,0}$ & $70(33 s^4-48 s^2+16)$ & $-385 s^4+500 s^2-136$\\
$1$ & $18(55 s^4-76 s^2+24)$  & $-3 \left(35 s^4-40 s^2+8\right)$       & $24(3 s^2-2)(5 s^2-4)$ & $0$ \\
$2$ & $-4(189 s^4-276 s^2+88)$ & $\Gamma_{2,2}$       & $-4(225 s^4-324 s^2+104)$ & $9(21 s^4-28 s^2+8)$ \\
$3$ & $-20(25 s^4-32 s^2+8)$ & $-\Lambda_{2,1}$       &    &  \\
\noalign{\smallskip}\hline
\end{tabular}
\end{table}

\section{Short period corrections} \label{a:shortperiod}

The short-period corrections of the new theory are the same as those reported in \citep{HautesserresLaraCeMDA2017}. For completeness, they are printed below. Thus,  $\{\xi,\mathcal{V}\}=-\frac{1}{2}J_2(R_\Earth/p)^2\Delta\xi$, where $\mathcal{V}$ denotes the generating function of the short-period elimination, and
\begin{align}
\Delta{r}=&\; p\left(1-\frac{3}{2}s^2-\frac{1}{2}s^2\cos2\theta\right), 
\\
\Delta\theta=&\; \Big[1-6c^2+(1-2c^2)\cos2\theta\Big]\sigma-\left[\frac{1}{4}-\frac{7}{4}c^2+(1-3c^2)\kappa\right]\sin2\theta,
\\
\Delta\nu=&\; c\left[(3+\cos2\theta)\sigma-\left(\frac{3}{2}+2\kappa\right) \sin2\theta\right],
\\
\Delta{R}=&\; \frac{\Theta}{r}\,(1+\kappa)s^2\sin2\theta,
\\
\Delta\Theta=& -\Theta\,s^2\left[\left(\frac{3}{2}+2\kappa\right)\cos2\theta+\sigma\sin2\theta\right],
\\
\Delta{N}=&\; 0,
\end{align}
where $\kappa$ and $\sigma$ denote the same functions of the polar variables previously defined in Eq.~(\ref{skpolar}).

\end{appendices}

\subsection*{Statements and Declarations}

\subsubsection*{Funding:} This work was supported by the Spanish State Research Agency and the European Regional Development Fund (Projects PID2021-123219OB-I00 and PID2024-161963OB-C22, AEI/ERDF, EU).

\subsubsection*{Conflict of Interest:} The author declares no conflict of interest.

\subsubsection*{Data Availability:} Data generated or analyzed in this study is contained within the article in the form of figures and tables.

\subsubsection*{Author Contribution:} Not applicable (single-authored paper).

\end{document}